\documentclass[11pt]{article}
\usepackage{amsmath}
\usepackage{amsfonts}
\usepackage{amssymb}
\usepackage{enumerate}
\usepackage{amsthm}
\usepackage{authblk}
\usepackage{graphicx}
\usepackage{verbatim}
\newtheorem{myproposition}{Proposition}[section]
\newtheorem{mytheorem}[myproposition]{Theorem}
\newtheorem{mylemma}[myproposition]{Lemma}

\newtheorem{mycorollary}[myproposition]{Corollary}

\newtheorem{mydefinition}[myproposition]{Definition}

\newtheorem{myproblem}[myproposition]{Problem}

\def\gr{\mathcal{G}}
\def\gchi{\chi^\Sigma_g}
\def\zet{\mathbb{Z}}

\makeatletter
\def\imod#1{\allowbreak\mkern10mu({\operator@font mod}\,\,#1)}
\makeatother
\begin {document}

\title{Group Sum Chromatic Number of Graphs\footnote{\copyright 2016. This manuscript version is made available under the CC-BY-NC-ND 4.0 license http://creativecommons.org/licenses/by-nc-nd/4.0/}}

\author[1,3]{Marcin Anholcer}
\author[2,3]{Sylwia Cichacz}

\affil[1]{\scriptsize{}Pozna\'n University of Economics and Business, Faculty of Informatics and Electronic Economy}
\affil[ ]{Al. Niepodleg{\l}o\'sci 10, 61-875 Pozna\'n, Poland, \textit{m.anholcer@ue.poznan.pl}}
\affil[ ]{}
\affil[2]{AGH University of Science and Technology, Faculty of Applied Mathematics}
\affil[ ]{Al. Mickiewicza 30, 30-059 Krak\'ow, Poland, \textit{cichacz@agh.edu.pl}}
\affil[ ]{}
\affil[3]{University of Primorska, Faculty of Mathematics, Natural Sciences and Information Technologies}
\affil[ ]{Glagolja\v{s}ka 8, SI-6000 Koper, Slovenia, \{marcin.anholcer, sylwia.cichacz-przenioslo\}@famnit.upr.si}

\maketitle

\begin{abstract}
We investigate the \textit{group sum chromatic number} ($\gchi(G)$) of graphs, i.e. the smallest value $s$ such that taking any Abelian group $\gr$ of order $s$, there exists a function $f:E(G)\rightarrow \gr$ such that the sums of edge labels properly colour the vertices. It is known that $\gchi(G)\in\{\chi(G),\chi(G)+1\}$ for any graph $G$ with no component of order less than $3$ and we characterize the graphs for which $\gchi(G)=\chi(G)$.
\end{abstract}


\section{Introduction}

A labelling of the edges of a graph $G$ is called a \textit{vertex-colouring} if it results in weighted degrees that properly colour the vertices. If we use the elements of $\{1,2,\dots,k\}$ to label the edges, such a labelling is called a \textit{vertex-colouring $k$-edge-labelling}.

The concept of colouring the vertices with the sums of edge labels was introduced for the first time by Karo\'nski, {\L}uczak and Thomason (\cite{ref_KarLucTho}). The authors posed the following question. Given a graph $G$ with no component of order less than $3$, what is the minimum $k$ such that there exists a vertex-colouring $k$-edge-labelling? As there are some analogies with ordinary proper graph colouring, we will call this minimum value of $k$ the \textit{sum chromatic number} and denote it by $\chi^\Sigma(G)$.

Karo\'nski, {\L}uczak and Thomason conjectured that $\chi^\Sigma(G)\leq 3$ for every graph $G$ with no component of order less than $3$. The first constant bound was proved by Addario-Berry et al. in \cite{ref_AddDalMcDReeTho} ($\chi^\Sigma(G)\leq 30$) and then improved by Addario-Berry et al. in \cite{ref_AddDalRee} ($\chi^\Sigma(G)\leq 16$), Wan and Yu in \cite{ref_WanYu} ($\chi^\Sigma(G)\leq 13$) and finally by Kalkowski, Karo\'nski and Pfender in \cite{ref_KalKarPfe2} ($\chi^\Sigma(G)\leq 5$).  

On the other hand, numerous authors studied various labelling problems, where elements of finite Abelian groups were used instead of integers to label  vertices or edges of graph. We give only a few examples here. Graham and Sloane in \cite{ref_GraSlo} studied harmonious graphs i.e. graphs for which there exists an injection $f:V(G)\rightarrow \mathbb{Z}_q$ that assigns to every edge $(x,y)\in E(G)$ a unique sum $f(x)+f(y)$ modulo $q$. Beals et al. (see \cite{ref_BeaGalHeaJun}) considered the concept of harmoniousness with respect to arbitrary Abelian groups. \.Zak in \cite{ref_Zak} generalized the problem and introduced a new parameter, the \textit{harmonious order of $G$}, defined as the smallest number $t$ such that there exists an injection $f:V(G)\rightarrow \mathbb{Z}_t$ (or surjection if $t<V(G)$) that produces distinct edge sums. Hovey in \cite{ref_Hov} considers the so-called $A-cordial$ labellings, where for a given Abelian group $A$ and a graph $G$ one wants to obtain a vertex labelling such that the classes of vertices labelled with one label are (almost) equinumerous and so are the classes of edges with the same sum. Cavenagh et al. (\cite{ref_CavComNel}) consider \textit{edge-magic total labellings} with finite Abelian groups, i.e. labellings of vertices and edges resulting in equal edge sums. Froncek in \cite{ref_Fro} defined the notion of group distance magic graphs, i.e. graphs allowing the bijective labelling of vertices with elements of an Abelian group resulting in constant sums of neighbour labels. Stanley in \cite{ref_Sta} studied the vertex-magic labellings of edges with the elements of an Abelian group $A$, i.e. labellings where the resulting weighted degrees are constant. Kaplan et al. in \cite{ref_KapLevRod} considered vertex-antimagic edge labellings, i.e. bijections $f:E(G)\rightarrow A$, where $A$ is a cyclic group, resulting in distinct weighted degrees of vertices. A reach survey of graph labellings, including group labellings, was published by Gallian \cite{ref_Gal}.

The problem considered in this paper is, in a sense, a connection between these two topics. Assume we are given an arbitrary graph $G$ that is $\chi(G)$-colourable and has no component isomorphic to $K_1$ or $K_2$. Assume $\gr$ is an Abelian group  of order $m\geq \chi(G)$ with the operation denoted by $+$ and neutral element $0$. For convenience we will write $ka$ to denote $a+a+\dots+a$ (where element $a$ appears $k$ times), $-a$ to denote the inverse of $a$ and we will use $a-b$ instead of $a+(-b)$. Moreover, the notation
$$\sum_{a\in S}{a}$$
will be used as a short form for $a_1+a_2+a_3+\dots$, where $a_1, a_2, a_3, \dots$ are all the elements of the set $S$.

Given an edge labelling $f:E(G)\rightarrow \gr$, we define the \textit{weighted degree} of a vertex $v$ as
$$
w(v)=\sum_{e\ni v}{f(e)}.
$$
We call $f$ a \textit{vertex-$\gr$-colouring} if the weighted degrees of adjacent vertices are distinct. The \textit{group sum chromatic number} of $G$, denoted by $\gchi(G)$, is the smallest integer $s$ such that for every Abelian group $\gr$ of order $s$ there exists a vertex-$\gr$-colouring labelling $f$ of $G$.

A partial solution to this problem was given by Karo\'nski, {\L}uczak and Thomason in \cite{ref_KarLucTho}. Namely, they proved that there always exists a \textit{vertex-$\gr$-colouring} labelling of a given graph $G$ if it is $|\gr|$-colourable and $|\gr|$ is odd, where $|\gr|$ stands for the order of $\gr$. However, for the sake of completeness we present our original proofs covering all the possible cases.

Let us define the following special family of graphs. Here and in the remainder of this article the expression \textit{odd color class} (\textit{even color class}) will mean \textit{color class of odd order} (\textit{colour class of even order}, respectively).

\begin{mydefinition}
A connected graph $G$ with order at least $3$ is called \emph{ugly} if it belongs to one of the following classes:
\begin{itemize}
\item
$\chi(G)=4k+2$ for some integer $k\geq 0$ and in each proper $\chi(G)$-colouring of $G$ all the colour classes of $G$ are odd,
\item
$\chi(G)=2^q$ for some integer $q\geq 2$ and in each proper $\chi(G)$-colouring of $G$ either exactly $2$ or exactly $2^q-2$ colour classes are odd.
\end{itemize} 
\end{mydefinition}

The main result of our paper is the following theorem, determining the value of $\gchi(G)$ for every graph with no component of order less than $3$.

\begin{mytheorem}\label{main_thm}
Let $G$ be an arbitrary graph with no component of order less than $3$. If $G$ does not contain any ugly component $C_1$ such that $\chi(C_1)=\chi(G)$ and, in the case when $\chi(G)=2^q$ for some integer $q\geq 2$, it has no component $C_2\cong K_{\chi(G)-2}$, then 
$$
\gchi(G)=\chi(G).
$$
\noindent{}Otherwise 
$$
\gchi(G)=\chi(G)+1.
$$
\end{mytheorem}

\section{Proof of Theorem \ref{main_thm}}

In order to properly colour all the vertices of $G$ we need at least $\chi(G)$ distinct elements of $\gr$. However it is not always enough, as shown by the following lemma. 

\begin{mylemma}\label{lemma_below}
If a connected graph $G$ is ugly, then
$$
\gchi(G)\geq\chi(G)+1.
$$
\end{mylemma}

\noindent\textbf{Proof.}
Assume first that $\chi(G)=4k+2$ and all the colour classes of $G$ are odd. Assume $\gchi(G)=4k+2$. Each Abelian group of order $4k+2$ is isomorphic to $\zet_2\times \gr^\star$ for some Abelian group $\gr^\star$ of order $2k+1$. We can express every element of $\gr$ as $(z,a)$, where $z\in \zet_2$ and $a\in\gr^\star$. In such a situation the elements of $2k+1$ of the colour classes would have weighted degrees equal to $(1,a_j)$ for distinct $a_j\in\gr^\star$ and the remaining $2k+1$ to $(0,a_j)$. This would imply that the sum of all the weighted degrees equals $(1,a)$ for some $a\in\gr^\star$. But on the other hand, for every edge label $(z,a_j)$ we have $2(z,a_j)=(0,2a_j)$, so $\sum_{v\in V(G)}w(v)=(0,b)$ for some $b\in\gr^\star$. This contradiction implies that $\gchi(G)\geq \chi(G)+1$.

Before we consider the second case, let us cite the following lemma (see e.g. \cite{ref_ComNelPal}, Lemma 8).

\begin{mylemma}\label{lemma_CNPsum0}
Let $\gr$ be an Abelian group.
\begin{enumerate}
\item
If $\gr$ has exactly one involution, say $i$, then
$$
\sum_{a\in \gr}{a}=\sum_{a\in \gr, 2a=0}{a}=i.
$$
\item
If $\gr$ has no involution, or more than one involution, then
$$
\sum_{a\in \gr}{a}=\sum_{a\in \gr, 2a=0}{a}=0.
$$
\end{enumerate}
\end{mylemma}  

\qed

Now we are ready to consider the second case. Assume that $\gr\cong \zet_2\times \zet_2\times\dots \times \zet_2$ ($q$ times, $q\geq 2$). It means that there are $2^q-1\geq 3$ involutions in $\gr$ and, according to Lemma \ref{lemma_CNPsum0}, their sum equals $0$. Assume that $\gchi(G)=2^q$. The vertices in every colour class must have weighted degrees equal to distinct elements of $\gr$. If there are two odd colour classes, let us assume that their vertices have weighted degrees $a$ and $b$, where $a\neq b$. The sum of the weighted degrees in any even class would be equal to $0$ and in the odd classes to $a$ and $b$. Thus $\sum_{v\in V(G)}w(v)=a+b$. Similarly, if there are exactly $2^q-2$ odd colour classes, then assume the weighted degrees of the vertices in two even classes are $a$ and $b$, where $a\neq b$. The sum of all the weighted degrees is equal to the sum of all the elements of $\gr$ except $a$ and $b$, that is $-a-b=a+b$. In both situations $\sum_{v\in V(G)}w(v)=0$. But the equality $a+b=0$ means $a=b$. The contradiction implies that $\gchi(G)\geq 2^q+1=\chi(G)+1$.

\qed

Given any two vertices $x_1$ and $x_2$ belonging to the same connected component of $G$, there exist walks from $x_1$ to $x_2$. Here by a walk  we mean a sequence of vertices, with possible repetitions, such that two consecutive vertices in this sequence are adjacent in $G$, together with the induced sequence of edges, also with possible repetitions. Some of the walks may consist of an even number of vertices (counting the repetitions). We are going to call them \textit{even walks}. The walks with an odd number of vertices are called \textit{odd walks}. We will always choose the shortest even or the shortest odd walk from $x_1$ to $x_2$.

We start with $0$ on all the edges of $G$. Then, in each step, we will choose $x_1$ and $x_2$ and add some labels to all the edges of the chosen walk from $x_1$ to $x_2$. To be more specific, we will add an element $a$ of the group to the labels of all the edges in odd position on the walk (starting from $x_1$) and $-a$ to the labels of all the edges in even position. It is possible that some labels will be modified more than once, as the walk does not need to be a path. We will denote such a situation by $\phi_e(x_1,x_2)=a$ if we label a shortest even walk and $\phi_o(x_1,x_2)=a$ if we label a shortest odd walk. Observe that putting $\phi_e(x_1,x_2)=a$ means adding $a$ to the weighted degrees of both $x_1$ and $x_2$, while $\phi_o(x_1,x_2)=a$ means adding $a$ to the weighted degree of $x_1$ and $-a$ to the weighted degree of $x_2$. In both cases the operation does not change the weighted degree of any other vertex of the walk.

Let us continue with the following lemma, determining the group sum chromatic number of connected bipartite graph. The upper bounds can be deduced from \cite[Theorem, p.152]{ref_KarLucTho} (ugly graphs) and \cite[Theorem 2.3]{ref_LuYuZha} (other graphs). However, as we already mentioned above, we present our proofs for the sake of completeness.

\begin{mylemma}\label{lemma_bipartite}
Let $G$ be a connected bipartite graph of order $n\geq 3$. Then
$$
\gchi(G)=\begin{cases}
3&\text{if $G$ is ugly,}\\
2&\text{otherwise}.
\end{cases}
$$
\end{mylemma}

\noindent\textbf{Proof.} We know from Lemma \ref{lemma_below} that the given values are lower bounds. Thus it suffices to present a proper group colouring. Let us denote the two colour classes by $V_1$ and $V_2$.

In the case when $G$ is not ugly, we have $\gr\cong \zet_2$ for every group $\gr$ of order $2$, so we can use $\zet_2$. At least one of the colour classes, say $V_1$ is even. We choose in any way $|V_1|/2$ disjoint ordered pairs of vertices $(x_1^j,x_2^j)$, $x_1^j,x_2^j\in V_1$ and put $\phi_o(x_1^j,x_2^j)=1$ for $j=1,\dots,|V_1|/2$. Finally we label all the unlabelled edges with $0$. This way we obtain
$$
w(x)=\begin{cases}
1&\text{if   } x\in V_1,\\
0&\text{otherwise}.
\end{cases}
$$
If $G$ is ugly (i.e. both colour classes are odd), then we can use $Z_3$, as for every Abelian group $\gr$ of order $3$ we have $\gr\cong Z_3$. Assume without loss of generality that $|V_1|\geq |V_2|$ and so $|V_1|\geq 3$. We choose three vertices $x_1,x_2,x_3\in V_1$. We choose in any way $(|V_1|-3)/2$ disjoint ordered pairs of vertices $(x_1^j,x_2^j)$, $x_1^j,x_2^j\in V_1\setminus\{x_1,x_2,x_3\}$ and put $\phi_o(x_1^j,x_2^j)=2$ for $j=1,\dots,(|V_1|-3)/2$. Then we put $\phi_o(x_1,x_2)=\phi_o(x_1,x_3)=1$. Finally we label all the unlabelled edges with $0$. This way we obtain
$$
w(x)=\begin{cases}
1\text{   or   }2&\text{if   } x\in V_1,\\
0&\text{otherwise}.
\end{cases}
$$

\qed

Before analyzing the case of non-bipartite graphs we need to prove the following technical lemma.

\begin{mylemma}\label{lemma_inv}
Let $\gr$ be an Abelian group with set of involutions $I^\star=\{i_1,\dots, i_{2^k-1}\}$, $k\geq 2$ and let $I=I^\star\cup \{0\}$. Then for any given $r$ such that $0\leq r \leq 2^k$, there exists a set $R\subseteq I$,  $|R|=r$, such that
$$
\sum_{i\in R}{i}=0
$$
if and only if $r\not\in \{2,2^k-2\}$.
\end{mylemma}

\noindent\textbf{Proof.} We can consider only the case when $r\leq 2^{k-1}$ because if we can find such a set $R$ for given $r$, then according to Lemma \ref{lemma_CNPsum0}, the set $I\setminus R$ will also sum up to $0$.

If $r=0$ then $R=\emptyset$ and if $r=1$, then $R=\{0\}$.

Let $r=2$. It is impossible to find $a,b\in I$, $a\neq b$ such that $a+b=0$, because the last condition implies $a=b$, a contradiction. Thus we cannot find a desired set $R$ neither for $r=2$, nor for $r=2^k-2$.

Now let $3\leq r\leq 2^{k-1}$. We know that $I\cong \zet_2\times \zet_2\times\dots\times \zet_2$, so each of its elements can be expressed as pairwise different binary sequence $i_r=(a_1^r,\dots,a_k^r)$. The addition in $I$ is represented by coordinate-wise addition in $\zet_2$. Assume these sequences are ordered lexicographically in increasing order. Observe that, for $q< 2^{k-1}$, we have $a_1^q=0$.

Let $R^\star=(i_1,\dots,i_{r-1})$ and let
$$
i^\star=\sum_{i\in R^\star}{i}.
$$

If $i^\star=0$, then $R=R^\star\cup\{0\}$ is the desired set. If $i^\star\neq 0$ and $i^\star=i_q$ for some $q>r-1$ (thus $i^\star\not\in R^\star$), then we put $R=R^\star\cup\{i^\star\}$. Finally, if $i^\star\neq 0$ and $i^\star=i_q$ for some $q\leq r-1$ ($i^\star\in R^\star$), then we choose any $i_p\in R^\star$, $p\neq q$ (such a $p$ exists as $|R^\star|=r-1\geq 2$) and define $R$ in the following way:
$$
R=R^\star\setminus \{i_p\}\cup \{i_p+i_{2^{k-1}},i_q+i_{2^{k-1}}\}.
$$
Indeed, we have
$$
\sum_{i\in R}{i}=\sum_{i\in R^\star}{i}-i_p+i_p+i_{2^{k-1}}+i_q+i_{2^{k-1}}=i^\star-i_p+i_p+i_{2^{k-1}}+i^\star+i_{2^{k-1}}=0,
$$
and $i_p+i_{2^{k-1}}=(1,a_2^p,\dots,a_k^p)\neq i_q+i_{2^{k-1}}=(1,a_2^q,\dots,a_k^q)$ are the only two elements of $R$ represented by the sequences with $1$ on the first position.

\qed

Before we proceed with the proof of Theorem \ref{main_thm} for non-bipartite graphs, let us introduce some notation.

We will denote the colour classes by $V_1,\dots,V_{\chi(G)}$, assuming that the odd classes precede the even ones in this sequence. We start with labelling all the edges with $0$. Now we are going to adjust some of the edge labels.

Observe that choosing any two vertices $x_i$ and $x_j$ of $G$ and any $a\in\gr$ , both operations $\phi_o(x_1,x_2)=a$ and $\phi_e(x_1,x_2)=a$ are allowed as $G$ is connected and contains at least one odd cycle.

Let us define two special ways of labelling the edges of the subsets of $V(G)$. Assume first that some colour class $V_j\subseteq V(G)$ is even. For any $a\in\gr$, we can increase the weighted degree of every vertex of $V_j$ by $a$ in the following way. If $a=0$, then we do not change the label of any edge. If $a\neq 0$, then we construct in any way $|V_j|/2$ disjoint pairs of vertices $(x_{i,1}^j,x_{i,2}^j)$, $i=1,\dots,|V_j|/2$ and for each such pair we put $\phi_e(x_{i,1}^j,x_{i,2}^j)=a$, increasing $w(x)$ by $a$ for every $x\in V_j$. In order to simplify the notation, we will denote such a labelling of $V_j$ by $\Phi(V_j)=a$. Observe that the operation $\Phi(V_j)=0$ is feasible also in the case when $V_j$ is odd.

Assume now that two colour classes $V_j,V_k\subseteq V(G)$ are odd. This time, for any couple $(-a,a)$ we can increase by $a$ the weighted degree of every vertex in $V_j$ and increase by $-a$ the weighted degree of every vertex in $V_k$. We do this by choosing arbitrary vertices $x_j\in V_j$ and $x_k\in V_k$ and putting $\phi_o(x_j,x_k)=a$, $\Phi(V_j\setminus\{x_j\})=a$ and $\Phi(V_k\setminus\{x_k\})=-a$. In order to simplify the notation, we will denote such an operation by $\Phi(V_j,V_k)=(a,-a)$.

The following result was proved by Karo{\'n}ski, {\L}uczak and Thomason (see \cite{ref_KarLucTho}, p.152).

\begin{mylemma}\label{lemma_oddGEn}
Let $G$ be a connected graph of order $n\geq 3$ and let $\gr$ be an arbitrary Abelian group of order $2k+1\geq \chi(G)$ for some integer $k$. Then there exists a vertex-$\gr$-colouring of $G$.
\end{mylemma}

\qed

Now we are ready to determine the value of $\gchi(G)$ for arbitrary connected non-bipartite graphs. In the remainder of the proof the number of odd classes will be denoted by $t$.

\begin{mylemma}\label{lemma_nobipartite}
Let $G$ be a connected non-bipartite graph of order $n\geq 3$. Then
$$
\gchi(G)=\begin{cases}
\chi(G)+1&\text{if $G$ is ugly,}\\
\chi(G)&\text{otherwise}.
\end{cases}
$$
\end{mylemma}

\noindent\textbf{Proof.} From Lemma \ref{lemma_below} it follows that it suffices to present the construction, proving this way the upper bound. 

If $\chi(G)$ is odd or $G$ is ugly, then $|\gr|$ is odd and we use Lemma \ref{lemma_oddGEn} to find the desired colouring.

If $\chi(G)$ is even and $G$ is not ugly, then $|\gr|=\chi(G)$ is even. If $t\leq 1$, then we assign distinct $a_j\in\gr$, $a_j\neq 0$ to the $|\gr|-1$ distinct even classes and put $\Phi(V_j)=a_j$, $j=2,\dots,|\gr|$. We also put $\Phi(V_1)=0$, no matter if $V_1$ is even or odd. If $t=2$, then we choose any $a_1\in\gr$, $a_1\neq 0$, $2a_1\neq 0$ (such an $a_1$ must exist as $G$ is not ugly and so $\gr\not\cong \zet_2\times\dots \times \zet_2$). Let us denote $a_2=-a_1$ and let $a_3,\dots,a_{|\gr|}$ be the remaining elements of $\gr$ ordered in any way. We put $\Phi(V_1,V_2)=(a_1,-a_1)$ and $\Phi(V_j)=a_j$ for $j=3,\dots,|\gr|$. Finally, if $t\geq 3$, then we have to distinguish four cases depending on the structure of $\gr$ and the exact value of $t$. Observe that $\gr\cong Z_{2^{p_1}}\times Z_{2^{p_2}}\times \dots\times Z_{2^{p_k}}\times \gr^\star$ for some positive integers $p_1$, $p_2$,\dots, $p_k$ and some Abelian group $\gr^\star$ where $|\gr^\star|=2p_0+1$ for some integer $p_0\geq 0$. There are exactly $2^k-1$ involutions $i_1,i_2,\dots,i_{2^k-1}$ in $\gr$. The remaining non-zero elements of $\gr$ form $(|\gr|-2^k)/2\geq 1$ disjoint pairs of the form $\{a_j,-a_j\}$.

\vspace{10pt}
\noindent{}\textit{Case $1$: }$k=1$.

\noindent{}We have exactly one involution $i_1$ in $\gr$. If $t<|\gr|$, then we can find $\lfloor t/2\rfloor$ disjoint couples of distinct non-zero elements the form $\{a_j,-a_j\}$, where $j=1,\dots, \lfloor t/2\rfloor$. We construct $\lfloor t/2\rfloor$ pairs of odd colour classes $(V_{2j-1},V_{2j})$, $1\leq j \leq\lfloor t/2\rfloor$ and for every such pair we put $\Phi(V_{2j-1},V_{2j})=(a_j,-a_j)$. Now, if $t$ is odd, then we put $\Phi(V_t)=0$. Finally, we assign to each of the even classes $V_j$, $j=t+1,\dots,|\gr|$ one of the remaining $|\gr|-t$ elements of $\gr$, $a_j$, and we put $\Phi(V_j)=a_j$.

If $t=|\gr|$, then $p_1\geq 2$ as $G$ is not ugly, and so $\gr$ has a subgroup $\gr^{\star\star}=\{0,a_1,2a_1=i_1,3a_1\}$. We choose any three vertices $x_1\in V_1$, $x_2\in V_2$ and $x_3\in V_3$ and put $\phi_e(x_1,x_2)=i_1$ and $\phi_e(x_1,x_3)=a_1$. Now, as $V_j\setminus\{x_j\}$ is even for every $j=1,2,3$, we put $\Phi(V_1\setminus\{x_1\})=3a_1$, $\Phi(V_2\setminus\{x_2\})=i_1$, $\Phi(V_3\setminus\{x_3\})=a_1$ and $\Phi(V_4)=0$. We order in pairs $(V_{2j-1},V_{2j})$, $j=3,\dots,|\gr|/2$, all the remaining $|\gr|-t$ odd classes. Then we arrange them with disjoint pairs of elements $(a_j,-a_j)$, $\{a_j,-a_j\}\subset \gr\setminus\gr^{\star\star}$ and we put $\Phi(V_{2j-1},V_{2j})=(a_j,-a_j)$.

\vspace{10pt}
\noindent{}\textit{Case $2$: }$k\geq 2, t\leq 2^k, t\neq 2^k-2$.

\noindent{}The number of involutions is $2^k-1\geq 3$. Using Lemma \ref{lemma_inv}, we choose $t$ elements $i_1,\dots,i_t\in \gr$ such that $2i_j=0$, $j=1,\dots,t$ and
$$
\sum_{j=1}^{t}{i_t}=0.
$$
We choose one vertex $x_j$ from every colour class $V_j$, $j=1,\dots,t$, and we put $\phi_o(x_1,x_j)=i_j$ for $j=2,\dots,t$. This way we obtain $w(x_j)=i_j$ for $j=1,\dots,t$. As in every $V_j$, $j=1,\dots,t$, the number of remaining vertices is now even, we put $\Phi(V_j\setminus\{x_j\})=i_j$, $j=1,\dots,t$. The remaining colour classes $V_j$ are even, so we can use the remaining $a_j\in\gr$ by putting $\Phi(V_j)=a_j$, $j=t+1,\dots,|\gr|$.

\vspace{10pt}
\noindent{}\textit{Case $3$: }$k\geq 2$, $t= 2^k-2$.

\noindent{}We choose any $a_1\in\gr$, $a_1\neq 0$, $2a_1\neq 0$ (such an $a_1$ must exist as $G$ is not ugly and so $\gr\not\cong \zet_2\times\dots \times \zet_2$). Let us denote $a_2=-a_1$. We put $\Phi(V_{t-1},V_{t})=(a_1,-a_1)$. Now we are left with $t-2=2^k-4\geq 4$ odd classes ($t\geq 3$ implies $k\geq 3$ in this case). Using Lemma \ref{lemma_inv}, we choose a set of $t-2$ elements $R=\{i_1,\dots,i_{t-2}\}\subseteq \gr$ such that $2i_j=0$, $j=1,\dots,t-2$ and
$$
\sum_{j=1}^{t-2}{i_t}=0.
$$
We choose an arbitrary vertex $x_j$ from every colour class $V_j$, $j=1,\dots,t-2$, and we put $\phi_o(x_1,x_j)=i_j$ for $j=2,\dots,t$. This way we obtain $w(x_j)=i_j$ for $j=1,\dots,t-2$. As in every $V_j$, $j=1,\dots,t-2$, the number of remaining vertices is now even, we put $\Phi(V_j\setminus\{x_j\})=i_j$, $j=1,\dots,t-2$. As the remaining colour classes $V_j$, $j=t+1,\dots,|\gr|$ are even, we can use the remaining $|\gr|-t$ elements $a_j\in\gr\setminus(R\cup\{a_1,a_2\})$ by putting $\Phi(V_j)=a_j$, $j=t+1,\dots,|\gr|$.

\vspace{10pt}
\noindent{}\textit{Case $4$: }$k\geq 2$, $t> 2^k$.

\noindent{}We can find $t_1=\lceil(t-2^k)/2\rceil$ disjoint couples of distinct non-zero elements of the form $\{a_j,-a_j\}$, where $j=1,\dots, t_1$. We construct $t_1$ pairs of odd colour classes $(V_{2j-1},V_{2j})$, $1\leq j \leq t_1$ and for every such pair we put $\Phi(V_{2j-1},V_{2j})=(a_j,-a_j)$.

Then we colour the remaining $t-2t_1$ odd classes as follows. Observe that either $t-2t_1=2^k-1$ or $t-2t_1=2^k$. Using Lemma \ref{lemma_inv}, we choose $t-2t_1$ elements $i_1,\dots,i_{t-2t_1}\in \gr$ such that $2i_j=0$, $j=1,\dots,{t-2t_1}$ and
$$
\sum_{j=1}^{t-2t_1}{i_j}=0.
$$
We choose one vertex $x_j$ from every colour class $V_j$, $j=2t_1+1,\dots,t$, and we put $\phi_o(x_{2t_1+1},x_j)=i_j$ for $j=2t_1+2,\dots,t$. This way we obtain $w(x_j)=i_j$ for $j=2t_1+1,\dots,t$. As in every $V_j$, $j=2t_1+1,\dots,t$, the number of remaining vertices is now even, we put $\Phi(V_j\setminus\{x_j\})=i_j$, $j=2t_1+1,\dots,t$.

The remaining $|\gr|-t$ colour classes $V_j$ are even, we can use the remaining $|\gr|-t$ elements $a_j\in\gr$ by putting $\Phi(V_j)=a_j$, $j=t+1,\dots,|\gr|$.

\qed

In order to finish the proof of the Theorem we need the following lemma.

\begin{mylemma}\label{lemma_extended}
Given a connected graph $G$ of order $n\geq 3$, it is possible to find a vertex-$\gr$-colouring of $G$ for every Abelian group $\gr$ of order $h>\gchi(G)$, except the case when $h=2^q$ for some integer $q\geq 2$ and $G=K_{h-2}$.
\end{mylemma}

\noindent\textbf{Proof.} Recall that $t$ denotes the number of odd colour classes of $G$. If $h$ is odd, then the existence of a vertex-$\gr$-colouring is guaranteed by Lemma \ref{lemma_oddGEn}. If $h$ is even, then we have to consider three cases.

\vspace{10pt}
\noindent{}\textit{Case $1$: }$|V(G)|> h$.

\noindent{}Using $h-\chi(G)$ additional colours, we can recolour some of the vertices in order to obtain a proper $h$-colouring of $G$. Moreover, we are able to control the parity of the new colour classes and obtain this way a ``good'' number of odd classes. Then we use Lemma \ref{lemma_nobipartite} in order to obtain the desired $\gr$-colouring.

\vspace{10pt}
\noindent{}\textit{Case $2$: }$|V(G)|\leq h$ and $G$ is bipartite. 

\noindent{}In this case, either one of the colour classes of $G$ is even, or both classes are odd. If there exists an even colour class, then we use the same method as in the proof of Lemma \ref{lemma_bipartite}, using one of the involutions of $\gr$ instead of $1$ (such an involution has to exist as $h$ is even). If both colour classes are odd and there is some element $a$ of order greater than $2$ in $\gr$, then we use the method from the proof of Lemma \ref{lemma_bipartite} substituting $a$ for $1$ and $-a$ for $2$. Finally, if both colour classes are odd and all the non-zero elements of $\gr$ are involutions, then we recolour the vertices of one of the colour classes of $G$ in order to obtain four odd colour classes $V_1,V_2,V_3,V_4$. Assume that the classes $V_2$, $V_3$ and $V_4$ are the new ones, so the bipartition of $G$ is $(V_1,V_2\cup V_3\cup V_4$). Observe that $h\geq 4$ as $\gchi(G)=3$ in this case. We apply Lemma \ref{lemma_inv} in order to properly colour the graph. To be more specific, we choose four distinct elements $i_j\in\gr$, $2i_j=0$, $j=1,2,3,4$ such that $\sum_{j=1}^{4}{i_j}=0$ and then we proceed similarly as in the proof of Lemma \ref{lemma_nobipartite} (the case when $3\leq t\leq 2^k$, $t\neq 2^k-2$): we choose any vertices $x_j\in V_j$, $j=1,2,3,4$ and we put $\phi_e(x_1,x_j)=i_j$, $j=2,3,4$ and $\Phi(V_j\setminus\{x_j\})=i_j$, $j=1,2,3,4$.

\vspace{10pt}
\noindent{}\textit{Case $3$: }$|V(G)|\leq h$ and $G$ is non-bipartite. 

\noindent{}In such a situation we have to distinguish two cases depending on the relation between $|V(G)|$ and $h$.

If $|V(G)|\neq h-2$ or $h\neq 2^q$ for any integer $q\geq 2$, then we recolour the vertices in order to obtain $t=|V(G)|$ odd (one-element) colour classes $V_j=\{x_j\}$, $j=1,\dots,t$. If there is only one involution in $\gr$, then we can find $\lfloor t/2\rfloor$ disjoint pairs $\{a_j,-a_j\}$, $a_j\neq 0$, $2a_j\neq0$ and we put $\phi_o(x_{2j-1},x_{2j})=a_j$, $j=1,\dots,\lfloor t/2\rfloor$. If the set of involutions $I^\star$ has $k\geq 2$ (i.e. $k\geq 3$) elements, then we can choose $t$ distinct elements $a_1,a_2,\dots,a_{t}\in\gr$ such that
$$
\sum_{j=1}^{t}{a_j}=0.
$$
To be more specific, if $(h-k-1)/2\geq \lfloor t/2\rfloor$, then we choose $\lfloor t/2\rfloor$ disjoint pairs $\{a_j,-a_j\}$, $a_j\neq 0$, $2a_j\neq0$, $j=1,\dots,\lfloor t/2\rfloor$. If $(h-k-1)/2< \lfloor t/2\rfloor$, then we choose all the $(h-k-1)/2$ disjoint pairs $\{a_j,-a_j\}\subset\gr$, $a_j\neq 0$, $2a_j\neq 0$, $j=1,\dots,(h-k-1)/2$ and we add to this set $t-(h-k-1)/2$ elements of $I^\star\cup\{0\}$ - their existence is guaranteed by Lemma \ref{lemma_inv}, if only $t-(h-k-1)/2\not\in \{2,k-1\}$. If $t-(h-k-1)/2\in \{2,k-1\}$, then we choose $(h-k-1)/2-1$ pairs $\{a_j,-a_j\}\subset\gr$ and $t-(h-k-1)/2+2$ elements of $I^\star\cup\{0\}$. The only situation when we cannot do that is when $(h-k-1)/2=0$, i.e. $\gr\cong \zet_2\times\dots\times \zet_2=(\zet_2)^q$ for some integer $q$ and $t\in\{2,h-2\}$ (i.e. $t=h-2$, as $V(G)\geq 3$).

Finally, if $|V(G)|= h-2$ and $\gr\cong \zet_2\times\dots\times \zet_2=(\zet_2)^q$ for some integer $q\geq 2$, then we consider two cases, depending on the structure of $G$.

If $G$ is not a complete graph, then $\chi(G)<|V(G)|$. That means that either the number of odd colour classes $t<|V(G)|$ satisfies $t\in\{0\}\cup\{3,\dots,h-3\}$, or there are exactly $t=2$ odd colour classes and we can recolour one of the even ones in order to obtain $t+2=4< h-2$ odd colour classes. In both situations we can apply the same methods as in the previous case in order to find a vertex-$\gr$-colouring labelling of $G$. Observe that we can exclude the case $t=1$, as it contradicts the parity of $|V(G)|=h-2$.

Assume finally that $G\cong K_{2^q-2}$. Suppose that there is some vertex-$\gr$-colouring labelling $f$ of $G$. Then the sum of all the weighted degrees would be $-a-b=a+b$ for some distinct elements $a,b\in\gr$, as the sum of all the elements of $\gr$ is $0$. On the other hand we have
$$
\sum_{x\in V(G)}{w(x)}=\sum_{e\in E(G)}{2f(e)}=0,
$$
so $a+b=0$ and $a=b$, a contradiction.

\qed

The main result of our paper follows from the above lemmas. If there is no ugly component $C$ such that $\chi(C)=\chi(G)$, then we know that $\gchi(G)\geq \chi(G)$. On the other hand, using Lemmas \ref{lemma_bipartite} and \ref{lemma_nobipartite} we are able to find a vertex-$\gr$-colouring edge labelling of every component $C_j$ of $G$, for any Abelian group $\gr$, $|\gr|=\gchi(C_j)\leq \chi(G)$. Thus by Lemma \ref{lemma_extended} we are able to find the desired $\gr$-colouring of $G$ for any Abelian group $\gr$, $|\gr|=\chi(G)$, except the case when there is some component $C_2\cong K_{\chi(G)-2}$ and $\chi(G)=2^q$ for some integer $q\geq 2$. But in such a case Lemma \ref{lemma_extended} guarantees that there exists a vertex-$\gr$-colouring edge labelling for every Abelian group $\gr$, $|\gr|\geq \chi(G)+1$.

If there is some ugly component $C$ such that $\chi(C)=\chi(G)$, then by Lemma \ref{lemma_below} we know that $\gchi(G)\geq \chi(G)+1$. On the other hand, using Lemmas \ref{lemma_bipartite} and \ref{lemma_nobipartite} we are able to find a vertex-$\gr$-colouring edge labelling of every component $C_j$ of $G$, for any Abelian group $\gr$, $|\gr|=\gchi(C_j)\leq \chi(G)+1$. Thus by Lemma \ref{lemma_extended} we are able to find the desired $\gr$-colouring of $G$ for any Abelian group $\gr$, $|\gr|=\chi(G)+1$ (note that $\chi(G)+1$ is odd).

\section{Final Remarks}

Using the argument similar to the one in the final part of the proof of Theorem \ref{main_thm}, we can formulate the following corollary.

\begin{mycorollary}
Let $G$ be any graph with no component of order less than $3$. Then for every $h\geq\gchi(G)$ there exists a vertex-$\gr$-colouring edge labelling of $G$ for every Abelian group $\gr$, $|\gr|=h$, except the case when $\gr\cong \zet_2\times\dots\times \zet_2\cong (\zet_2)^q$ for some integer $q\geq 1$ and there exists a component $C\cong K_{h-2}$ of $G$.
\end{mycorollary}

In the introduction we mentioned that the problem considered here is similar to the problem of finding $\chi^\Sigma(G)$, i.e. the smallest $k$ such that there exists a vertex-colouring $k$-edge-labelling. Also a total version of this problem has been defined, where vertex labelling is allowed and the weighted degree is calculated as

$$
w(v)=f(v)+\sum_{e\ni v}{f(e)}
$$

\noindent{}(see e.g. \cite{ref_Kal}). In the case of the problem presented in this paper, such a modification does not make sense. Trivially, for every graph $G$, any group of order $\chi(G)$ would be enough as one can label the vertices of different colour classes with distinct elements of $\gr$ and then put $0$ on all the edges.
 
In the proof of Theorem~\ref{main_thm} we often use the fact that we are allowed to use $0$ on edges. Thus a natural question is the following.

\begin{myproblem}
Let $G$ be a simple graph with no component of order less than $3$. For any Abelian group $\gr$, let $\gr^\star = \gr \setminus \{0\}$. Determine the  \textit{non-zero group sum chromatic number} (${\gchi}^\star(G)$) of $G$, i.e. the smallest value $s$ such that taking any Abelian group $\gr$ of order $s$, there exists a function $f:E(G)\rightarrow \gr^\star$ such that the resulting weighted degrees properly colour the vertices.
\end{myproblem}

In this case the variant of the problem with labels on vertices allowed does not have to be trivial.

All the elements of $\gr$ can be obtained as some combination of not necessarily all of its elements, in particular of its generators. The question is, how many elements of $\gr$ we have to use in order to obtain a vertex-$\gr$-colouring edge labelling.

\begin{myproblem}
Assume that for a given simple graph $G$ with no component of order less than $3$ there exists a vertex-$\gr$-colouring edge labelling for every group $\gr$ of order $s$. What is the minimum number $k=k(G,s)$ such that for every group $\gr$ of order $s$ there is a subset $S\subseteq\gr$, $|S|\leq k$ such that there exists a vertex-$\gr$-colouring edge labelling $f:E(G)\rightarrow S$?
\end{myproblem}

So far we considered only finite Abelian groups. So, the next question seems to be natural, as some generalization of the problem of the ordinary vertex-colouring edge-labelling.

\begin{myproblem}
Let $G$ be a simple graph with no component of order less than $3$. Determine the smallest value $k$ such that for any infinite Abelian group $\gr$ there exists a subset $S\subseteq\gr$, $S\leq k$ such that there exists a vertex-$\gr$-colouring labelling $f:E(G)\rightarrow S$.  
\end{myproblem}

Observe that $\gchi(G)=\chi(G)+1$ does not mean that for each group $\gr$ of order $|\gr|=\chi(G)$ there is no vertex-$\gr$-colouring of $G$. It would be interesting to characterize such groups of order $|\gr|=\chi(G)$, which still allow a vertex-$\gr$-colouring of $G$ with $\gchi(G)=\chi(G)+1$. In particular, one may consider the following open problem.

\begin{myproblem}
Assume that $\gchi(G)=\chi(G)+1=2^q+1$ for a graph $G$, where $q\geq 2$ is an integer. Determine all groups of order $|\gr|=2^q$ for which a vertex-$\gr$-colouring exists.  
\end{myproblem}

\subsection*{Acknowledgments}
We are grateful to the leaders, researchers and administrative staff of the Faculty of Mathematics, Natural Sciences and Information Technologies of the University of Primorska for creating the perfect work conditions during our post-doc fellowship and inspiring discussions. We are also grateful to Prof. Dalibor Froncek for inspiration. Finally, we are indebted to the anonymous referees for their remarks that allowed us to improve this paper.


\begin{thebibliography}{99}

\bibitem{ref_AddDalMcDReeTho}
L. Addario-Berry, K. Dalal, C. McDiarmid, B.A. Reed, A. Thomason, Vertex-colouring edge-weightings, Combinatorica 27 (2007) 1--12.

\bibitem{ref_AddDalRee}
L. Addario-Berry, K. Dalal, B.A. Reed, Degree constrainted subgraphs, Discrete Applied Mathematics 156 (2008) 1168--1174.

\bibitem{ref_BeaGalHeaJun}
R. Beals, J. Gallian, P. Headley, D. Jungreis, Harmonious groups, Journal of Combinatorial Theory Series A 56 (1991) 223--238.

\bibitem{ref_CavComNel}
N. Cavenagh, D. Combe, A.M. Nelson, Edge-magic group labellings of countable graphs, Electronic Journal of Combinatorics 13 (2006) {\#}R92, 19 pp.

\bibitem{ref_ComNelPal}
D. Combe, A.M. Nelson, W.D. Palmer, Magic labellings of graphs over finite Abelian groups, Australasian Journal of Combinatorics 29 (2004) 259--271.

\bibitem{ref_Fro}
D. Froncek, Group distance magic labeling of of Cartesian products of cycles, Australasian Journal of Combinatorics 55 (2013) 167--174.

\bibitem{ref_Gal}
J.A. Gallian, A dynamic survey of graph labeling, Electronic Journal of Combinatorics, DS6, http://www.combinatorics.org/Surveys/.

\bibitem{ref_GraSlo}
R.L. Graham, N.J.A. Sloane, On additive bases and harmonious graphs, SIAM Journal on Algebraic and Discrete Methods, 1 (1980) 382--404.

\bibitem{ref_Hov}
M. Hovey, $A$-cordial graphs, Discrete Mathematics 93 (1991) 183--194.

\bibitem{ref_Kal}
M. Kalkowski, A note on 1,2-conjecture, Electronic Journal of Combinatorics (submitted for publication).

\bibitem{ref_KalKarPfe2}
M. Kalkowski, M. Karo\'nski, F. Pfender, Vertex-coloring edge-weightings: Towards the 1-2-3-conjecture, Journal of Combinatorial Theory Series B 100 (2010) 347--349.

\bibitem{ref_KapLevRod}
G. Kaplan, A. Lev, Y. Roditty, On zero-sum partitions and anti-magic trees, Discrete Mathematics 309 (2009) 2010--2014.

\bibitem{ref_KarLucTho}
M. Karo\'nski, T. {\L}uczak, A. Thomason, Edge weights and vertex colours, Journal of Combinatorial Theory Series B 91 (2004) 151--157.

\bibitem{ref_LuYuZha}
H. Lu, Q. Yu, C.-Q. Zhang, Vertex-coloring 2-edge-weighting of graphs, European Journal of Combinatorics 32 (2011) 21--27.

\bibitem{ref_Sta}
R.P. Stanley, Linear homogeneous Diophantine equations and magic labelings of graphs, Duke Mathematical Journal 40 (1973) 607--632.

\bibitem{ref_WanYu}
T. Wang, Q. Yu, On vertex-coloring 13-edge-weighting, Frontiers of Mathematics in China 3 (2008) 581--587.

\bibitem{ref_Zak}
A. \.Zak, Harmonious orders of graphs, Discrete Mathematics 309 (2009) 6055--6064.

\end{thebibliography}
\end {document}